\documentstyle[12pt]{amsart}

\newsymbol\boxtimes 1202

\newcommand{\nc}{\newcommand}

\nc{\ad}{\operatorname{ad}}
\nc{\Aut}{\operatorname{Aut}}
\nc{\Boxtimes}{\fbox{$\times$}} 
\nc{\blt}{\bullet}
\nc{\bSt}{\mbox{\bf{St}}}
\nc{\card}{\operatorname{card}}
\nc{\Cch}{\check{C}}
\nc{\cd}{\operatorname{cd}}
\nc{\Ch}{\operatorname{Ch}}
\nc{\chara}{\operatorname{char}}
\nc{\CHom}{\cal{H}om}
\nc{\Coker}{\operatorname{Coker}}
\nc{\codim}{\operatorname{codim}}
\nc{\Cone}{\operatorname{Cone}}
\nc{\cSgn}{\cal{S}gn}
\nc{\depth}{\operatorname{depth}} 
\nc{\dirlim}{\underset{\rightarrow}{\operatorname{lim}}}
\nc{\Div}{{\operatorname{Div}}}
\nc{\dlog}{{\operatorname{dlog}}}
\nc{\dotbox}{\overset{\bullet}{\boxtimes}}
\nc{\dotimes}{\overset{\bullet}{\otimes}}
\nc{\emp}{\emptyset}
\nc{\Ext}{\operatorname{Ext}}
\nc{\Fac}{\cal{F}ac}
\nc{\FM}{{\cal{FM}}}
\nc{\Fun}{\operatorname{F}}
\nc{\FS}{{{\cal{{FS}}}}}
\nc{\gr}{{\operatorname{gr}}}
\nc{\Hom}{\operatorname{Hom}}
\nc{\hgt}{\operatorname{ht}}
\nc{\hsl}{\widehat{{\frak{sl}}}}
\nc{\Id}{\operatorname{Id}}
\nc{\id}{\operatorname{id}}
\nc{\Ima}{\operatorname{Im}}
\nc{\ind}{\operatorname{ind}}
\nc{\Ind}{\operatorname{Ind}}
\nc{\infi}{\operatorname{inf}}
\nc{\infh}{\frac{\infty}{2}}
\nc{\invlim}{\underset{\leftarrow}{\operatorname{lim}}}
\nc{\Ker}{\operatorname{Ker}}
\nc{\Locsys}{\cal{L}ocsys}
\nc{\Mod}{\operatorname{Mod}}
\nc{\modul}{\operatorname{mod}}
\nc{\Ob}{\operatorname{Ob}}
\nc{\opp}{\operatorname{opp}}
\nc{\Or}{\cal{O}r}
\nc{\Ord}{\cal{O}rd}
\nc{\Part}{\cal{P}art}
\nc{\PGL}{\operatorname{PGL}}
\nc{\Qui}{{\operatorname{Qui}}}
\nc{\sgn}{\operatorname{sgn}}
\nc{\Sh}{\cal{S}h}
\nc{\sll}{{\frak{sl}}}
\nc{\Spe}{{\operatorname{Sp}}}
\nc{\supr}{\operatorname{sup}}
\nc{\Supp}{\operatorname{Supp}}
\nc{\supp}{\operatorname{supp}}
\nc{\tFM}{{\widetilde{\cal{FM}}}}
\nc{\tFS}{\widetilde{\cal{FS}}}
\nc{\Tor}{\operatorname{Tor}}
\nc{\totimes}{\tilde{\otimes}}
\nc{\tr}{\operatorname{tr}}
\nc{\Vect}{\cal{V}ect}
\nc{\wt}{\widetilde}

\nc{\bo}{\mbox{\bf{0}}}
\nc{\One}{\mbox{\bf{1}}}
\nc{\one}{\mbox{\bf{1}}}
  
\nc{\BA}{{\Bbb{{A}}}}
\nc{\ba}{\mbox{\bf{a}}}
\nc{\baJ}{\bar{J}}
\nc{\BAO}{\overset{\circ}{\Bbb A}}
\nc{\BB}{\Bbb B}
\nc{\bB}{\mbox{\bf{B}}}
\nc{\bc}{\mbox{\bf{c}}}
\nc{\BC}{{\Bbb C}}
\nc{\bCC}{\bar{\cal{C}}}
\nc{\bCW}{{\overline{\cal{W}}}}
\nc{\bD}{\bar{D}}
\nc{\bd}{\mbox{\bf{d}}}
\nc{\BE}{\overline{E}}
\nc{\BF}{\overline{F}}
\nc{\bF}{\mbox{\bf{F}}}
\nc{\bof}{\mbox{\bf{f}}}
\nc{\bL}{\mbox{\bf{L}}}
\nc{\blambda}{\bar{\lambda}}
\nc{\bM}{\mbox{\bf{M}}}
\nc{\bmu}{\vec{\mu}}
\nc{\BN}{\Bbb{N}}
\nc{\bnu}{\vec{\nu}}
\nc{\BP}{{\Bbb P}}
\nc{\bP}{\mbox{\bf{P}}}
\nc{\BPO}{\overset{\circ}{\Bbb{P}}}
\nc{\BQ}{\Bbb{Q}}
\nc{\bq}{\mbox{\bf{q}}}
\nc{\BR}{\Bbb R}
\nc{\br}{\mbox{\bf{r}}}
\nc{\breta}{\bar{\eta}}
\nc{\bs}{\mbox{\bf{s}}}
\nc{\bt}{\mbox{\bf{t}}}
\nc{\bU}{\mbox{\bf{U}}}
\nc{\bu}{\mbox{\bf{u}}} 
\nc{\BUpsilon}{\bar{\Upsilon}}
\nc{\bw}{\mbox{\bf{w}}}
\nc{\bx}{\mbox{\bf{x}}}
\nc{\BZ}{{\Bbb Z}}
\nc{\bz}{\mbox{\bf{z}}}

\nc{\CA}{{\cal{A}}}
\nc{\CB}{{\cal{B}}}
\nc{\CC}{{\cal{C}}}
\nc{\CD}{{\cal{D}}}
\nc{\CE}{{\cal{E}}}
\nc{\CF}{{\cal{F}}}
\nc{\CH}{{\cal{H}}}
\nc{\CI}{{\cal{I}}}
\nc{\CJ}{{\cal{J}}}
\nc{\CK}{{\cal{K}}}
\nc{\CL}{{\cal{L}}}
\nc{\CM}{{\cal{M}}}
\nc{\CN}{{\cal{N}}}
\nc{\CO}{{\cal{O}}}
\nc{\CP}{{\cal{P}}}
\nc{\CQ}{{\cal{Q}}}
\nc{\CR}{{\cal{R}}}
\nc{\CS}{{\cal{S}}}
\nc{\CT}{{\cal{T}}}
\nc{\CU}{{\cal{U}}}
\nc{\CV}{{\cal{V}}}
\nc{\CW}{{\cal{W}}}
\nc{\CX}{{\cal{X}}}
\nc{\CY}{{\cal{Y}}}
\nc{\CZ}{{\cal{Z}}}

\nc{\dd}{\operatorname{d}}
\nc{\DO}{\overset{\circ}{D}}
\nc{\dpar}{\partial}

\nc{\fA}{\frak{A}}
\nc{\fE}{\frak{E}}
\nc{\fF}{\frak{F}}
\nc{\ff}{\frak{f}}
\nc{\fg}{\frak{g}}
\nc{\fh}{{\frak{h}}}
\nc{\fl}{\frak{l}}
\nc{\fn}{{\frak{n}}}
\nc{\fp}{\frak{p}}
\nc{\fu}{\frak{u}}

\nc{\HO}{\overset{\circ}{H}}
\nc{\hfg}{\hat{\frak{g}}}
\nc{\hfn}{\hat{\frak{n}}}
\nc{\hL}{\hat{L}}

\nc{\jo}{\overset{\circ}{j}}

\nc{\phid}{\overset{\bullet}{\phi}}

\nc{\tBP}{\tilde{\Bbb{P}}}
\nc{\tC}{\tilde{C}}
\nc{\tc}{\tilde{c}}
\nc{\tCA}{\tilde{\cal{A}}}
\nc{\tCC}{\tilde{\cal{C}}}
\nc{\tCI}{\tilde{\cal{I}}}
\nc{\tCO}{\tilde{\cal{O}}}
\nc{\tCP}{\tilde{\cal{P}}}
\nc{\tCT}{\tilde{\cal{T}}}
\nc{\tD}{\tilde{D}}
\nc{\tDelta}{\tilde{\Delta}}
\nc{\tE}{\tilde{E}}
\nc{\tF}{\tilde{F}}
\nc{\tfF}{\tilde{\frak{F}}}
\nc{\tff}{\tilde{\frak{f}}}
\nc{\tfu}{\tilde{\frak{u}}}
\nc{\tJ}{\tilde{J}}
\nc{\tj}{\tilde{j}}
\nc{\tK}{\tilde{K}}
\nc{\tL}{\tilde{L}}
\nc{\tM}{\tilde{M}}
\nc{\tP}{\tilde{P}}
\nc{\tPhi}{\tilde{\Phi}}
\nc{\TPO}{\overset{\circ}{T\Bbb{P}}}
\nc{\tR}{\tilde{R}}
\nc{\tS}{\tilde{S}}
\nc{\tT}{\tilde{T}}
\nc{\ttau}{\tilde{\tau}}
\nc{\ttheta}{\tilde{\theta}}
\nc{\tU}{\tilde{U}}
\nc{\tUpsilon}{\tilde{\Upsilon}}
\nc{\ty}{\tilde{y}}
\nc{\tY}{\tilde{Y}}
\nc{\txi}{\tilde{\xi}}

\nc{\UD}{\overset{\bullet}{U}}
\nc{\UO}{\overset{\circ}{U}}
\nc{\UU}{\operatorname{U}}

\nc{\valpha}{\vec{\alpha}}
\nc{\vbeta}{\vec{\beta}}
\nc{\vc}{\vec{c}}
\nc{\vD}{\vec{D}}
\nc{\vd}{\vec{d}} 
\nc{\vgamma}{\vec{\gamma}}
\nc{\vK}{\vec{K}}
\nc{\vlambda}{\vec{\lambda}}
\nc{\vmu}{\vec{\mu}}
\nc{\vnu}{\vec{\nu}}
\nc{\vo}{\vec{0}}
\nc{\vu}{\vec{u}}
\nc{\vx}{\vec{x}}

\nc{\XO}{\overset{\circ}{X}}

\nc{\nen}{\newenvironment}
\nc{\ol}{\overline}
\nc{\ul}{\underline}
\nc{\ra}{\rightarrow}
\nc{\lra}{\longrightarrow}
\nc{\Lra}{\Longrightarrow}
\nc{\lla}{\longleftarrow}
\nc{\Llra}{\Longleftrightarrow}
\nc{\hra}{\hookrightarrow}
\nc{\iso}{\overset{\sim}{\lra}}
\nc{\rlh}{\rightleftharpoons}

\nc{\Thm}[1]{Theorem~\ref{#1}}
\nc{\Prop}[1]{Proposition~\ref{#1}}
\nc{\Lem}[1]{Lemma~\ref{#1}}
\nc{\Cor}[1]{Corollary~\ref{#1}}
\nc{\Conj}[1]{Conjecture~\ref{#1}}
\nc{\Claim}[1]{Claim~\ref{#1}}
\nc{\Defn}[1]{Definition~\ref{#1}}
\nc{\Exa}[1]{Example~\ref{#1}}
\nc{\Rem}[1]{Remark~\ref{#1}}
\nc{\Note}[1]{Note~\ref{#1}}

\nen{thm}[1]{\label{#1}{\bf Theorem.\ } \em}{}
\nen{prop}[1]{\label{#1}{\bf Proposition.\ } \em}{}
\nen{lem}[1]{\label{#1}{\bf Lemma.\ } \em}{}
\nen{cor}[1]{\label{#1}{\bf Corollary.\ } \em}{}
\nen{conj}[1]{\label{#1}{\bf Conjecture.\ } \em}{}
\nen{claim}[1]{\label{#1}{\bf Claim.\ } \em}{}

\nen{defn}[1]{\label{#1}{\bf Definition.\ } }{}
\nen{exa}[1]{\label{#1}{\bf Example.\ } }{}

\nen{rem}[1]{\label{#1}{\em Remark.\ } }{}
\nen{note}[1]{\label{#1}{\em Note.\ } }{}
\nen{exer}[1]{\label{#1}{\em Exercise.\ } }{}

\begin{document}

\title[Rational differential forms on the line]{Rational differential forms on the line 
and singular vectors in Verma modules over $\hsl_2$}

\author[Vadim Schechtman and Alexander Varchenko]
{Vadim Schechtman$^{*}$
and Alexander Varchenko$^{\diamond,1}$}

\thanks{${}^1\,$ Supported in part by NSF grant DMS-1362924 and Simons Foundation grant \#336826.}

\maketitle

\begin{center}
{\it ${}^*$ Institut de Math\'ematiques de Toulouse\,--\,  Universit\'e Paul Sabatier\\ 118 Route de Narbonne,
31062 Toulouse, France \/}

\medskip
{\it $^\diamond$Department of Mathematics, University of North Carolina
at Chapel Hill\\ Chapel Hill, NC 27599-3250, USA\/}
\end{center}

\begin{abstract}
We construct a monomorphism of  the De Rham complex of scalar multivalued
meromorphic forms on the projective line, holomorphic
on the complement to a finite set of points, to 
the chain complex of the Lie algebra  of $\sll_2$-valued
algebraic functions on the same complement with coefficients in
a tensor product of contragradient Verma modules over the affine Lie
algebra $\hsl_2$.  We  show that the existence of singular vectors
 in the Verma modules (the Malikov-Feigin-Fuchs singular vectors)
is reflected in the new relations between the cohomology classes
of logarithmic differential forms.

\end{abstract}

\section{Introduction} 

\subsection{}

We consider two complexes. The first 
is the De Rham complex of scalar multivalued
meromorphic forms on the projective line that are holomorphic
on the complement to a finite set of points. The second
is the chain complex of the Lie algebra  of $\sll_2$-valued
algebraic functions on the same complement with coefficients in
a tensor product of contragradient Verma modules over the affine Lie
algebra $\hsl_2$. We construct a monomorphism of the first
complex to the second and show that the existence of singular vectors in the Verma modules
is reflected in the new relations between the cohomology classes
of logarithmic differential forms.

This construction has two  motivations.

The first motivation was
to generalize the principal construction of \cite{sv}.
In \cite{sv},  we identified the tensor products
of contragradient Verma modules over a semisimple
Lie algebra and the spaces of the top degree
{ logarithmic}  differential forms over certain configuration spaces.
We also identified the logarithmic parts of the De Rham
complexes over the configuration spaces with some  standard Lie algebra chain
complexes having coefficients in these tensor products, cf. in
\cite{ks} a  $\CD$-module explanation of this correspondence.

The second idea was that the appearance of { singular
vectors} in Verma modules over affine Lie algebras is
reflected in the {\em new relations}
between the cohomology classes of logarithmic differential forms;
moreover, in some sense this correspondence should be one-to-one. This
was proved in an important particular case in \cite{fsv}; in \cite{stv}
a one-to-one correspondence was established "on the level of parameters".
In the present work we construct (see Section \ref{res sing}) this correspondence for another
non-trivial class of singular vectors, namely for (a part of)
{ Malikov-Feigin-Fuchs} (MFF) ones, cf. \cite{mff}. (The first examples,
were worked out during the preparation of \cite{fsv}.)
It turns out that MFF vectors (having quite complicated form) admit
a very simple definition using certain {\em limiting procedure}, see
Theorem \ref{limit}.
      
Our paper is related to the recent paper \cite{AFO} by M.\,Aganagic, E.\,Frenkel,  A.\,Okounkov devoted to 
quantum q-Langlands correspondence.  In Section 6 of \cite{AFO} the authors discuss how  conformal blocks of a WZW model  are related to conformal blocks of the \lq\lq{}dual\rq\rq{} $\mathcal W$-algebra. If the conformal blocks are defined by one-dimensional integrals the problem is reduced to comparing multivalued
meromorphic forms on the projective line in terms of representation theory of $\hsl_2$.

\medskip
This paper had been written in the middle of 90s and prepared for publication in the fall
 of 2015 while the second author visited the MPI in Bonn. The second author thanks
MPI for hospitality.
The authors thank E.\,Mukhin for interest in this paper, E.\,Frenkel for interesting discussions, and the anonymous referee who helped to improve the exposition.

\section{The De Rham complex of a hypergeometric function}

\subsection{}
\label{def dr}  
Let $z_1, \ldots, z_n$ be pairwise distinct complex numbers, $z_{n+1}=\infty$, 
$U=\BC-\{z_1,\ldots,z_n\}$, 
\begin{equation}
\ell=\prod_{1\leq i<j\leq n}\ (z_i-z_j)^{M^iM^j/2\kappa}\ 
\prod_{i=1}^n\ (t-z_i)^{-M^i/\kappa}\ ,
\end{equation}
where $t$ is a coordinate on $\BC$ and $M^1,\ldots,M^n,\kappa$ are complex parameters.
 The function $\ell$ 
is a multivalued  holomorphic function on $U$ with 
singularities at $z_1,\ldots,z_n$ and infinity. 

The function $\ell$ defines a hypergeometric function 
of $z_1,\ldots,z_n$ by the formula
\begin{equation}
I(z_1,\dots,z_n)=\int_{\gamma}\ \ell\ dt.
\end{equation}
Here $\gamma$ is a suitable cycle on $U$, for example 
a path connecting two points $z_i$, $z_j$. 

Consider the {\em twisted De Rham complex} associated 
with $\ell$:
\begin{equation}
\label{dr com}
0\lra\Omega^0\overset{\tilde d}{\lra}\Omega^1\lra 0 .
\end{equation}
Here $\Omega^p$ is the space of rational differential forms 
on $\BC$ regular on $U$. The differential $\tilde d$ is given 
by the formula
\begin{equation}
\label{nabla}
\tilde d=d_{DR}+\alpha\wedge\cdot\,,
\end{equation}
where $d_{DR}$ is the De Rham differential and the second 
summand is the left exterior multiplication by the form
\begin{equation}
\label{alpha}
\alpha=-\frac{1}{\kappa}\sum_{i=1}^n\ \frac{M^i}{t-z_i}dt.
\end{equation}
Formula (\ref{nabla}) is motivated by the computation 
$$
d_{DR}(\ell\omega)=\ell d_{DR}\omega+d_{DR}\ell\wedge\omega
=\ell(d_{DR}\omega+\alpha\wedge\omega).
$$
The complex  
$\Omega^{\bullet}$ is the complex of global {\em algebraic} sections of 
the De Rham complex of $(\CO^{an}_U,\nabla)$ where 
$\nabla=d_{DR}+\alpha\wedge\cdot$ is the integrable 
connection  
on the sheaf $\CO^{an}_U$ of holomorphic functions on $U$.  

If $\CS$ is the locally constant sheaf 
of horizontal sections then the cohomology $H^{\bullet}(U;\CS)$ is equal to 
$H^{\bullet}(\Omega^{\bullet})$. 

If the monodromy of $\ell$ is non-trivial, then 
\begin{equation}
H^0(\Omega^{\bullet})=0,\ \dim\ H^1(\Omega^{\bullet})=n-1.
\end{equation}

\subsection{}
\label{log forms}  
The simplest elements of $\Omega^1$ are {\em logarithmic forms}: 
\begin{equation}
\label{log}
\omega_i=M^i\frac{d(t-z_i)}{t-z_i}\ , \qquad i=1,\, \ldots, \,n.
\end{equation}
They are cohomologically dependent:
\begin{equation}
\omega_1+\ldots+\omega_n=-\kappa d(1).
\end{equation}
For generic $M^1,\ldots, M^n,\kappa$ the forms $\omega_1,\ldots,\omega_n$ 
generate the space $H^1$ and the relation $\sum\ \omega_i\sim 0$ is the 
only one, \cite{stv}.

\subsection{Resonances} 
\label{reson} For special {\em resonance} values of parameters 
the forms $\omega_1,\ldots,\omega_n$ span a proper subspace of $H^1$. 

Here are the resonance conditions.

(a) $M^i=-a\kappa$ where $a=0,1,2,\ldots;\ i=1,\ldots,n$. 

(b) $M^{n+1}=-2+a\kappa$ where $a=1,2,\ldots$. Here 
$M^{n+1}:=M^1+\ldots+M^n-2$. 

Each resonance condition implies a new cohomological relation between
the forms $\omega_1,\ldots,\omega_n$.  

\subsection{Example} 
\label{e1}
\label{ex 1} If $M^{n+1}=-2+\kappa$, then $\sum_{i=1}^n z_i\omega_i
\sim 0$. 

\subsection{Example}
\label{e2} 
\label{ex 2} If $M^{n+1}=-2+2\kappa$, then 
$$
\sum_{i=1}^n\ z_i^2\omega_i-\frac{1}{\kappa}\left(\sum_{j=1}^n\ 
z_jM^j\right)\left(\sum_{i=1}^n\ z_i\omega_i\right)\sim 0.
$$

It turns out that there is a direct connection between the
relations among the forms $\omega_1,\ldots,\omega_n$ 
and singular vectors in Verma modules over the 
affine Lie algebra $\hsl_2$. An explanation of this relation is 
the subject of this work. 

\subsection{} 
\label{restr} For any $i=1,\ldots,n$, define a number 
$a_i$ by the formula $a_i:=-M^i/\kappa$,  if $-M^i/\kappa$ is a 
non-negative integer,  and set $a_i:=\infty$ otherwise. 

Set $a_{n+1}:=(M^1+\ldots+M^n)/\kappa$ if the right-hand side 
is a positive integer and $a_{n+1}:=\infty$ otherwise. 

The number $\card \{i\in\{1,\ldots,n+1\}|\ a_i<\infty\}$ 
will be called the {\em number of resonances}. 

Introduce the {\em restricted De Rham complex}
\begin{equation}
0\lra\Omega_R^0\lra\Omega_R^1\lra 0
\end{equation}
as the subcomplex of (\ref{dr com}) where $\Omega^p_R
\subset\Omega^p$ is the subspace of the forms $\omega$ such that 

({\bf{r}}) for any $i=1,\ldots,n+1$, the degree of 
the pole of 
$\omega$ at the point $t=z_i$ is not greater than $a_i$ 
if $a_i<\infty$. 

\subsection{Lemma} 
\label{res forms} ${}$

{\em{\em (a)} $\Omega^{\bullet}_R$ is a 
subcomplex of $\Omega^{\bullet}$.

{\em (b)} The forms $\omega_1,\ldots,\omega_n$ belong to 
$\Omega^1_R$ and generate the space $H^1(\Omega^{\bullet}_R)
$. 

{\em (c)} The natural homomorphism $H^1(\Omega^{\bullet}_R)
\lra H^1(\Omega^{\bullet})$ is a monomorphism. 

{\em (d)} The codimension of the subspace $H^1(\Omega^{\bullet}_R)\subset H^1(\Omega^{\bullet})$ 
is equal to the number of resonances. 

{\em (e)} The forms $d(t-z_i)/(t-z_i)^{-a_i-1}$ for 
resonance points $t=z_i$ and the form $t^{a_{n+1}-1}dt$, if 
$t=\infty$ is a resonance point, give a basis of the 
space $H^1(\Omega^{\bullet})/H^1(\Omega^{\bullet}_R)$.} 

\subsection{} It is convenient to use the following 
basis in $\Omega^{\bullet}$. 

(a) An {\em elementary function} is a function 
$(t-z_i)^{-b}\ (b\in\BZ_{>0})$ or $t^b\ (b\in\BZ_
{\geq 0})$. 

(b) An {\em elementary differential form} is a form 
$d(t-z_i)/(t-z_i)^b\ (b\in\BZ_{>0})$ or $t^bdt\ (b\in\BZ_
{\geq 0})$. 

We have two basic formulas:  
\begin{equation}
\label{bff}
\kappa \,\tilde d((t-z_i)^{-b})=-(M^i+b\kappa)d(t-z_i)/(t-z_i)
^{b+1}+
\end{equation}
$$\sum_{k=1}^b\ \sum_{j\neq i}\ M^j/(z_j-z_i)^k\cdot 
d(t-z_i)/(t-z_i)^{b+1-k}-\sum_{j\neq i}\ M^j/(z_j-z_i)^b
\cdot d(t-z_j)/(t-z_j)
$$
and 
\begin{equation}
\label{bfi}
\kappa \,\tilde d(t^b)=\left(b\kappa -\sum_{j=1}^n\ M^j\right)t^{b-1}dt-
\end{equation}
$$
\sum_{k=1}^{b-1}\ \sum_{j=1}^n\ M^jz_j^kt^{b-1-k}dt-
\sum_{j=1}^n\ M^jz_j^b d(t-z_j)/(t-z_j) .
$$
If the resonance condition \ref{reson}(a) is satisfied 
and $b=a$, then the first term in the right-hand side of (\ref{bff}) 
disappears. Similarly, if the condition \ref{reson}(b) 
is satisfied and $b=a$, then the first term in the right-hand side 
of (\ref{bfi}) disappears. 

In the next sections we will give an interpretation 
for the elementary functions (forms) and formulas 
(\ref{bff}), (\ref{bfi}) in terms of 
$\hsl_2$-representations. 

Lemma \ref{res forms} follows easily from (\ref{bff}) and 
(\ref{bfi}).

\section{The Gauss-Manin connection}
\label{GMC}
In this section we show the important fact that the subbundle with fiber 
$H^1(\Omega^{\bullet}_R(z_1,\dots,z_n))$ $\subset H^1(\Omega^{\bullet}(z_1,\dots,z_n))$ 
is invariant with respect to the Gauss-Manin connection on the bundle with fiver
$H^1(\Omega^{\bullet}(z_1,\dots,z_n))$,
which moves the points $z_1,\dots,z_n$.
By Lemma \ref{res forms} the classes of logarithmic forms $\omega_1,\dots,\omega_n$ generate 
$H^1(\Omega^{\bullet}_R(z_1,\dots,z_n))$ for every distinct $z_1,\dots,z_n$. Our goal will be to describe
the relations between these classes in terms of $z_1,\dots,z_n$.

\subsection{} 
\label{fiber} When the points $z_1,\ldots,z_n$ are moving, the cohomology 
groups $H^{\bullet}(U;\CS)$ (as well as the dual homology groups 
$H_{\bullet}(U;\CS^*)$) form a vector bundle with a flat {\em Gauss-Manin} 
connection.

Set $\BC^{[n]}:=\BC^n-\bigcup_{i<j}\ \{z\in\BC^n\ |\ z_i=z_j\}$; 
$\BC^{[n+1]}:=\BC^{n+1}-\bigcup_{i<j}\ \{(z,t)\in\BC^{n+1}\ |\ z_i=z_j\}-
\bigcup_{i=1}^n\ \{(z,t)\in\BC^{n+1}\ |\ t=z_i\}$. 
Let $\psi: \BC^{[n+1]}\lra\BC^{[n]},\ (z,t)\mapsto z$, be the projection. 
This projection is a locally trivial bundle with  fiber 
$U(z)=\{t\in\BC\ |\ t\neq z_i,\ i=1,\ldots,n\}$.

Define an integrable connection $\nabla=d_{DR}+\beta\wedge\cdot$ on the 
sheaf $\CO^{an}$ of holomorphic functions on $\BC^{[n+1]}$,
where $\beta\wedge\cdot$ is the left multiplication by 
the form 
\begin{equation}
\beta=-\sum_{i=1}^n\ \frac{M^i}{\kappa}\ \frac{d(t-z_i)}
{t-z_i}+\sum_{i<j}\ \frac{M^iM^j}{2\kappa}\ 
\frac{d(z_i-z_j)}{z_i-z_j},
\end{equation}
cf. (\ref{alpha}). Let $\CS$ be the locally constant 
sheaf of horizontal sections. 
The fiber bundle $\psi$ together with the local system 
$\CS$ on $\BC^{[n+1]}$ defines a vector bundle 
$R^1\psi$ on $\BC^{[n]}$ with the fiber $H^1(U(z);\CS|_{U(z)})$ 
over $z\in\BC^{[n]}$. We have an isomorphism
$$
H^1(U(z);\CS|_{U(z)})=H^1(\Omega^{\bullet}(z)) ,
$$
where $\Omega^{\bullet}(z)$ is the twisted De Rham complex 
of the fiber defined in \ref{def dr}. 

This vector bundle has a canonical {\em Gauss-Manin connection} 
\begin{equation}
\nabla=\sum\ \nabla_{z_i}dz_i ,
\end{equation}
which can be defined as follows. Let $A\subset\BC^{[n+1]}$ 
be a Zariski open set, $\Omega^1(A\times\BC)$ the space 
of rational differential forms on $A\times\BC$ whose 
poles are on the hyperplanes $t=z_i$. A form $\omega
\in\Omega^1(A\times\BC)$ defines a section $[\omega]$ 
of $R^1\psi$ over $A$ with the value 
$[\omega|_{U(z)}]\in H^1(\Omega^{\bullet}(z))$ at $z$. 
The form $\eta:=d_{DR}\omega+\beta\wedge\omega$ can be 
written as 
\begin{equation}
\label{eta}
\eta=\sum_{i<j}\ \eta_{ij}(t,z)dz_i\wedge dz_j+
\sum_i\ \eta_i(t,z)dz_i\wedge dt ,
\end{equation}
where $\eta_{ij}, \eta_i$ are functions. By definition,
$$
 \nabla_{z_i}[\omega]:=[\eta_i dt].
$$
Elementary differential forms $d(t-z_i)/(t-z_i)^b$ and 
$t^bdt$ generate $\Omega^1(A\times\BC)$ over the ring 
of rational functions on $\BC^{[n]}$ regular on $A$. 
Hence, to compute the Gauss-Manin connection it is 
sufficient to compute (\ref{eta}) for elementary 
differential forms. 

The following two formulas give a description of the 
Gauss-Manin connection:
\begin{equation}
\kappa\beta\wedge\frac{d(t-z_i)}{(t-z_i)^b}=
\end{equation}
$$
\sum_{j<k;\ i\not\in\{j,k\}}\ \frac{M^jM^k}{2}\ 
\frac{d(z_j-z_k)}{z_j-z_k}\wedge  
\frac{d(t-z_i)}{(t-z_i)^b}
+\sum_{j\neq i}\ 
\frac{M^j(M^i-2)}{2}\ \frac{d(z_j-z_i)}{(z_j-z_i)}\wedge  
\frac{d(t-z_i)}{(t-z_i)^b}+
$$
$$
\sum_{j\neq i}\ M^j\frac{d(z_j-z_i)}{(z_j-z_i)^b}\wedge
\frac{d(t-z_j)}{t-z_j}-\sum_{j\neq i}\ \sum_{m=1}^{b-1}\ 
M^j\frac{d(z_i-z_j)}{(z_j-z_i)^{b-m+1}}\wedge
\frac{d(t-z_i)}{(t-z_i)^m},
$$
\begin{equation}
\kappa\beta\wedge t^bdt=
\sum_{j<k}\ \frac{M^jM^k}{2}\ \frac{d(z_j-z_k)}{(z_j-z_k)}
\wedge t^bdt+
\end{equation}
$$
\sum_{a=0}^{b-1}\ \sum_{i=1}^n\ M^iz_i^{b-a+1}dz_i\wedge 
t^adt+\sum_{i=1}^n\ M^iz_i^b dz_i\wedge\frac{d(t-z_i)}
{t-z_i}.
$$
The first formula has an important special case. For 
the logarithmic forms defined in (\ref{log}), we have 
\begin{equation}
\label{cov log}
\kappa\beta\wedge\omega_i=
\sum_{j<k;\ i\not\in\{j,k\}}\ \frac{M^jM^k}{2}\ 
\frac{d(z_j-z_k)}{z_j-z_k}\wedge\omega_i+
\end{equation}
$$
\sum_{j\neq i}\ \frac{M^j(M^i-2)}{2}\ 
\frac{d(z_j-z_i)}{z_j-z_i}\wedge\omega_i+
\sum_{j\neq i}\ M^i\frac{d(z_j-z_i)}{z_j-z_i}\wedge\omega_j.
$$
\subsection{Corollary} {\em Consider the subbundle 
$R^1\psi_R\subset R^1\psi$ with  fiber $H^1(\Omega^{\bullet}_R(z))
\subset H^1(\Omega^{\bullet}(z))$. This subbundle is 
invariant under the Gauss-Manin connection.} 

Here $\Omega^{\bullet}_R(z)$ is the restricited De Rham 
complex introduced in \ref{restr}. 

In fact, by Lemma \ref{res forms}(b), the logarithmic 
forms generate the fibers of $R^1\psi_R$ and by 
(\ref{cov log}), their covariant derivatives are 
expressed through logarithmic forms.

\section{Representations of $\hsl_2$}

\subsection{} Let $\sll_2$ be the Lie algebra of complex $2\times 2$-matrices 
with the zero trace; let $e,f,h$ be its standard generators, subject 
to the relations
$$
[e,f]=h,\quad [h,e]=2e,\quad  [h,f]=-2f.
$$
Let $\hsl_2$ be the corresponding affine Lie algebra 
$\hsl_2=\sll_2[T,T^{-1}]\oplus\BC c$ with the bracket 
$$
[aT^i,bT^j]=[a,b]T^{i+j}+i\langle a,b\rangle\delta_{i+j,0} c ,
$$
where $c$ is the central element, $\langle a,b\rangle:=\tr (ab)$. 

Set $f_1=f,\ e_1=e,\ h_1=h,\ f_2=eT^{-1},\ e_2=fT,\ h_2=c-h$. These 
elements are the standard Chevalley generators defining $\hsl_2$ as the 
Kac-Moody algebra corresponding to the Cartan matrix
$$
\left(\begin{array}{cc}2&-2\\-2&2\end{array}\right) .
$$
\subsection{Remark} 
\label{aut} The algebra $\hsl_2$ has an 
automorphism $\pi$ sending $\ c, eT^i, fT^i, hT^i$ to 
$\ c$, $fT^i$, $eT^i$, $-hT^i$ respectively. 

\subsection{} Fix $k\in\BC$, the value of the central charge. 
We assume that the action of $c$ on all our representations 
is the multiplication by $k$. 

For $M\in\BC$, let $V(M,k-M)$ be the {\em Verma module} 
over $\hsl_2$ generated by a vector $v$ subject to 
the relations $e_1v=e_2v=0,\ h_1v=Mv,\ h_2v=(k-M)v$. 

Let $\UU\hfn_-$ be the enveloping algebra of the Lie 
subalgebra $\hfn_-\subset\hsl_2$ generated by $f_1, f_2$. 
The map $\UU\hfn_-\lra V(M,k-M)$, $F\mapsto Fv$, is an 
isomorphism of $\UU\hfn_-$-modules. The space $V(M,k-M)$ 
has a $\Gamma=\BZ^2_{\geq 0}$-grading: 
a vector $f_{i_1}\cdot\ldots\cdot f_{i_p}v$\ with $i_j\in\{1,2\}$ 
has grading $(p_1,p_2)$, where $p_i$ is the number 
of $i$'s in the sequence $i_1,\ldots,i_p$. 

For $\gamma\in\Gamma$,  denote by $V(M,k-M)_{\gamma}\subset V(M,k-M)$ the corresponding $\gamma$-homogeneous component.
A homogeneous vector $\omega$ in 
$V(M,k-M)$, non-proportional to $v$, is 
called a {\em singular vector} if $e_1\omega=e_2\omega=0$. The Verma module
 $V(M,k-M)$ is reducible if and only if it contains 
a singular vector. For generic $M$, $V(M,k-M)$ is 
irreducible.

\subsection{Reducibility conditions}
\label{reduc} (See Kac-Kazhdan  
\cite{kk}). Set $\kappa:=k+2$. The Verma module 
$V(M,k-M)$ is 
reducible if and only if at least one of the equations (a)-(c)  
below is satisfied. 

(a) $M=l-1-(a-1)\kappa$. 

(b) $M=-l-1+a\kappa$.

(c) $\kappa=0$. 

Here $l,a=1,2,3, \ldots$ . If $(M,\kappa)$ satisfies 
exactly one of the conditions (a), (b), then the module 
$V(M,k-M)$ contains a unique proper submodule and 
this submodule is generated by a singular vector 
of degree $(la,l(a-1))$ for  condition (a) and of degree
$(l(a-1),la)$ for  condition (b). 

The singular vectors are highly nontrivial and are given 
by the following theorem.

\subsection{Theorem} (Malikov-Feigin-Fuchs, \cite{mff}) 
{\em For any positive integers $a, l$ and  $\kappa\in 
\BC$, the monomial 

$\operatorname{(a)}$ \phantom{aaaaa} $F_{12}(l,a,\kappa)=f_1^{l+(a-1)\kappa}f_2^{l+(a-2)\kappa}
f_1^{l+(a-3)\kappa}\cdot\ldots\cdot f_2^{l-(a-2)\kappa}
f_1^{l-(a-1)\kappa}$ 

lies in $\UU\hfn_-$. If  $M=l-1-(a-1)\kappa$, then  $F_{12}(l,a,\kappa)v$ is a singular 
vector of  $V(M,k-M)$ of degree $(la,l(a-1))$.
Similarly, the monomial

$\operatorname{(b)}$ $F_{21}(l,a,\kappa)=f_2^{l+(a-1)\kappa}f_1^{l+(a-2)\kappa}
f_2^{l+(a-3)\kappa}\cdot\ldots\cdot f_1^{l-(a-2)\kappa}
f_2^{l-(a-1)\kappa}$ 

lies in $\UU\hfn_-$.  If  $M=-l-1+a\kappa$, then  $F_{21}(l,a,\kappa)v$ is a singular 
vector of $V(M,k-M)$ of degree $(l(a-1),la)$. }

The explanation of the meaning of complex powers in these formulas see in \cite{mff}.

\subsection{Examples} 
\label{ex 3} {\bf 1.} $M=-2+\kappa$, 
$F_{21}(1,1,\kappa)v=f_2v=\frac{e}{T}v$. 

{\bf 2.} $M=-2+2\kappa$, 
$$
F_{21}(1,2,\kappa)v=f_2^{1+\kappa}
f_1f_2^{1-\kappa}v=f\left(\frac{e}{T}\right)^2v+(1+\kappa)\frac{h}{T}\ 
\frac{e}{T}v-(1+\kappa)\kappa\frac{e}{T^2}v.
$$
\subsection{Claim} The reducibility conditions for $l=1$ 
(the dimension of the complex line is $1$) correspond 
to the resonance conditions for $H^1(\Omega^{\bullet}_R)
\subset H^1(\Omega^{\bullet})$, and the singular vectors 
correspond to the relations between the forms 
$\omega_1,\ldots,\omega_n$, cf. examples \ref{ex 1}, 
\ref{ex 2} and \ref{ex 3}. 

We will make this statement precise in the next section. 

\subsection{} 
\label{contra} The maximal proper submodule of a Verma 
module $V$ coincides with the kernel of the {\em Shapovalov 
form}, which is the unique symmetric bilinear 
form $S(\cdot,\cdot)$ on $V$ characterized by the 
conditions $S(v,v)=1,\ S(f_ix,y)=S(x,e_iy)$ for all $i=1,2;\ 
x,y\in V$. 

One can regard $S$ as a map $S:\ V\lra V^*$ where 
$V^*:=\oplus_{\gamma\in\Gamma}\ V_{\gamma}^*$ and
$V^*_{\gamma}$ being the dual space to $V_{\gamma}$. 
There is a unique  $\hsl_2$-module structure on $V^*$ such 
that $\langle f_i\phi,x\rangle=\langle \phi,e_ix\rangle;\ 
\langle e_i\phi,x\rangle=\langle\phi,f_ix\rangle $, where $\phi\in V^*,\ x\in V,\ 
i=1,2 $. We call this $\hsl_2$-module $V^*$  the 
{\em contragradient dual} of $V$. The map $S$ 
is a morphism of $\hsl_2$-modules. The quotient 
$L:=V/\Ker\ S$ is irreducible.

\section{The main homomorphism}

{\em Conformal block construction.}  

\subsection{} 
\label{conf block} In this section let  $z=(z_1,\ldots,z_n,z_{n+1}=\infty)$ be pairwise distinct points 
of the complex projective line $\BP^1$. Fix local coordinates $t-z_1,\ldots, 
t-z_n,1/t$ at these points. Set $U(z):=\BP^1-\{z_1,\ldots,z_{n+1}\}$.  Notice that this $U(z)$ is the same $U(z)$ as in Section \ref{GMC}.

Let $\sll_2(U(z))$ be the Lie algebra of $\sll_2$-valued rational functions 
on $\BP^1$ regular on $U(z)$, with the pointwise bracket. 
Let $W_1,\ldots,W_{n+1}$ be representations of $\hsl_2$. The algebra 
$\sll_2(U(z))$ acts on the space $W_1\otimes\ldots\otimes W_{n+1}$: 
\begin{equation}
\label{action} 
a(t)\cdot(w_1\otimes\ldots\otimes w_{n+1}) \mapsto [a(t+z_1)] w_1
\otimes w_2\otimes\ldots\otimes w_{n+1}+ \dots +
\end{equation}
$$
 w_1\otimes\ldots\otimes w_{n-1}\otimes [a(t+z_n)] 
w_n\otimes 
 w_{n+1}+ w_1\otimes\ldots\otimes w_n\otimes\pi([a(1/t)]) w_{n+1} ,
$$
where $[b(t)]$ denotes the Laurent expansion of a function $b(t)$ at $t=0$
and the letter  
$\pi$ denotes the automorphism of $\hsl_2$ introduced in \ref{aut}.
We assume that all representations $W_i$ 
have the following finiteness property: 

({\bf Fin}) given $ w\in W_i$ and $ a\in\sll_2$, we have $aT^j\cdot w=0$ for all $j \gg 0$, 

so that the action of Laurent power series is well defined.  
The action 
of $c$ adds up to zero due to the residue formula. 
 Thus, we have the multiplication map
\begin{equation}
\label{mult}
\mu(z):\ \sll_2(U(z))\otimes (\otimes_{i=1}^{n+1}\ W_i)\lra 
\otimes_{i=1}^{n+1}\ W_i\ .
\end{equation}
The space
\begin{equation}
\label{sp conf}
(\otimes_{i=1}^{n+1}\ W_i)_{\sll_2(U(z))}:=\Coker\ \mu(z)
\end{equation} 
is called the {\em space of conformal blocks} at $z$. 

{\em The Knizhnik-Zamolodchikov connection}, see  \cite{kz}, \cite{f}. 

\subsection{} Let $\{X^a\},\ a=1,2,3$, be an orthonormal 
basis of $\sll_2$. Set
$$
L_{-1}:=\frac{1}{\kappa}\sum_{i=0}^{\infty}\ 
\sum_{a=1}^3\ (X^aT^{-i-1})(X^aT^i).
$$
It is a well defined operator on a representation 
satisfying the property (Fin) above. 

We have 
\begin{equation}
[L_{-1},XT^i]=-iXT^{i-1}
\end{equation}
for any $X\in\sll_2$ and any $i$. 

\subsection{} The following notation will be used: the action of 
an element $X$ of an algebra on the $i$-th factor 
of a tensor product of modules will be denoted by 
$X^{(i)}$. 

Recall the notations of \ref{fiber}. For an open subspace $A\subset\BC^{[n]}$, 
set $U_A:=(A\times\BP^1)\cap\BC^{[n+1]}\subset\BC^{[n]}\times\BP^1$. Denote 
by $\sll_2(U_A)$ the Lie algebra of algebraic $\sll_2$-valued functions 
on $U_A$. 

Let $W_1,\ldots,W_{n+1}$ be representations of $\hsl_2$ satisfying (Fin). 
Consider the trivial vector bundle $\CW_A:=A\times (W_1\otimes\ldots\otimes 
W_{n+1})\lra A$.   
The Lie algebra $\sll_2(U_A)$ acts on its holomorphic sections by formula 
(\ref{action}). 

Consider the flat connection on the bundle $\CW_A$: $\nabla=\sum_{i=1}^n\ 
\nabla_{z_i}\ dz_i$, 
\begin{equation}
\label{kz conn}
\nabla_{z_i} G(z)=\dpar_{z_i}G(z)+L_{-1}^{(i)}G(z),
\end{equation} 
where $G(z)\in\Gamma(A;\CW_A)$. 

\subsection{Lemma} (Cf. \cite{f}) {\em For any $X\in\sll_2(U_A),\ 
G\in\Gamma(A;\CW_A)$, we have 
$\nabla_{z_i}(XG)=(\nabla_{z_i}X)G+X(\nabla_{z_i}G)$.} 

This is proved by a direct computation. 

\subsection{} Let Im $\mu\subset\CW_{\BC^{[n]}}$ denote the subspace
whose intersection with the fiber at  a point $z$ is equal to 
the image of $\mu (z)$. According to the lemma this subspace is 
invariant with respect to the connection. Consider the quotient
bundle $\bCW$ over $\BC^{[n]}$ with fiber Coker $\mu (z)$.
This bundle is called {\em the bundle of conformal blocks}.

\subsection{Corollary} {\em The connection defined in
(\ref{kz conn}) induces a  connection on the bundle of
conformal blocks.}

The induced integrable connection on $\bCW$ is called the {\em Knizhnik-Zamolodchikov} 
(KZ) connection. The {\em KZ equation} on $G$ is the horizontality 
condition, $\nabla G=0$. 

{\em The main construction.} 

\subsection{} Let $M^1,\ldots,M^n, k$ be complex numbers,  $k\neq -2$. Set 
$M^{n+1}:=M^1+\ldots+M^n-2$. Let $V_i$ denote the Verma module $V(M^i,k-M^i)$ 
and $V_i^*$ the contragradient dual, cf. \ref{contra}. 

According to \ref{conf block}, the Lie algebra $\sll_2(U(z))$ acts on 
$\otimes_{i=1}^{n+1}\ V_i^*$, so that we can consider the standard 
chain complex $C_{\bullet}(\sll_2(U(z));\otimes_{i=1}^{n+1}\ V_i^*)$. Its 
right end looks as in (\ref{mult}):
\begin{equation}
\label{chains}
C_{\bullet}(\sll_2(U(z));\otimes_{i=1}^n\ V_i^*):\ \ldots\lra
\sll_2(U(z))\otimes\ (\otimes_{i=1}^{n+1}\ V_i^*)\overset{d}{\lra}
\otimes_{i=1}^{n+1}\ V_i^*\lra 0, 
\end{equation}
with $d=\mu(z)$, where $\mu(z)$ is defined in  \ref{conf block}. We assign to the last term degree $0$ and agree that $d$ has 
degree $1$, so that the whole complex sits in the nonpositive area.  

On the other hand, consider the shifted De Rham complex (\ref{dr com}) 
corresponding to $\kappa=k+2$: 
\begin{equation}
\Omega^{\bullet}(U(z))[1]:\ 0\lra\Omega^0(U(z))\lra\Omega^1(U(z))\lra 0.
\end{equation}
Here the shift $[1]$ means simply that we assign to $\Omega^j$ degree $j-1$. 

In the rest of the section we  construct a monomorphism of complexes 
$\Omega^{\bullet}(U(z))[1]\hra C_{\bullet}(\sll_2(U(z));\otimes_{i=1}^{n+1}\ 
V_i^*)$.

\subsection{}
\label{BASIS}
 First we need a basis in the complex 
(\ref{chains}). Let $\gamma=(p_1,p_2)\in\BZ^2_{\geq 0}$ 
and $p_1>p_2$.  
 We fix the following bases of  homogeneous components $V_{\gamma}$ of all Verma modules $V$: 
\begin{equation}
\frac{f}{T^{i_1}}\cdot\ldots\cdot\frac{f}{T^{i_a}}
\frac{h}{T^{j_1}}\cdot\ldots\cdot\frac{h}{T^{j_b}}
\frac{e}{T^{l_1}}\cdot\ldots\cdot\frac{e}{T^{l_c}}v
\end{equation}
where
\begin{equation}
\label{indices} 
0\leq i_a\leq i_{a-1}\leq\ldots\leq i_1,\ 
1\leq j_b\leq j_{b-1}\leq\ldots\leq j_1,\ 
1\leq l_c\leq l_{c-1}\leq\ldots\leq l_1;\ 
\end{equation}
$$
\sum_{s=1}^a\,i_s\,+\,\sum_{s=1}^b\,j_s\,+\,\sum_{s=1}^c\, l_s\, +\,a
\,-\,c\,=\,p_1, \qquad 
\sum_{s=1}^a\,i_s\,+\,\sum_{s=1}^b\,j_s\,+\,\sum_{s=1}^c\, l_s\,=\,p_2.
$$
For $p_1<p_2$, we fix a basis of the form 
\begin{equation}
\frac{e}{T^{l_1}}\cdot\ldots\cdot\frac{e}{T^{l_c}} 
\frac{h}{T^{j_1}}\cdot\ldots\cdot\frac{h}{T^{j_b}}
\frac{f}{T^{i_1}}\cdot\ldots\cdot\frac{f}{T^{i_a}}v,
\end{equation}
with the indices satisfying (\ref{indices}). These are bases 
by the Poincar\'{e}-Birkhoff-Witt theorem.

Notice that the elements $X/T^i$ and $X/T^j\ (X\in\sll_2)$ 
commute. 

We fix the bases in the contragradient Verma modules $V^*$
which are dual to the bases distinguished above in the Verma modules.

If $\{v_i\}$ is a basis in $V$, then we denote the
dual basis by $\{(v_i)^*\}$.

\subsection{} Define a map 
$$
\eta^1:\ \Omega^1(U(z))\lra \otimes_{i=1}^{n+1}\ V_i^*
$$
by the formulas
\begin{equation}
\label{eta 1f}
\frac{d(t-z_m)}{(t-z_m)^{b+1}}\mapsto 
-\kappa v_1^*\otimes\ldots\otimes \left(\frac{f}{T^b}v_m\right)^*
\otimes\ldots\otimes v_{n+1}^*,\ 
\end{equation}
\begin{equation}
\label{eta 1i}
t^bdt\mapsto \kappa v_1^*\otimes\ldots\otimes v_n^*
\otimes \left(\frac{e}{T^{b+1}}v_{n+1}\right)^*,
\end{equation}
for $b\geq 0$. 

Define a map
$$
\eta^0:\ \Omega^0(U(z))\lra \sll_2(U(z))\otimes 
(\otimes_{i=1}^{n+1}\ V_i^*)
$$
by the formulas
\begin{equation}
\label{eta 0f}
\frac{1}{(t-z_m)^b} \, \mapsto \, \frac{f}{(t-z_m)^b}\otimes 
v_1^*\otimes\ldots\otimes v_{n+1}^*-
\end{equation}
$$
\sum_{l=1}^b\ \Big[\frac{e}{(t-z_m)^l}\otimes v_1^*\otimes 
\ldots\otimes 2\sum_{i+j=b-l,\ i\geq j\geq 0}\ 
\left(\frac{f}{T^i}\frac{f}{T^j}v_m\right)^*\otimes\ldots\otimes 
v_{n+1}^*+
$$
$$
\frac{h}{(t-z_m)^l}\otimes v_1^*\otimes
\ldots\otimes \left(\frac{f}{T^{b-l}}v_m\right)^*\otimes
\ldots\otimes v_{n+1}^* \Big],
$$
for $b\geq 1$;
\begin{equation}
\label{eta 0i}
1 \, \mapsto \, f\otimes v_1^*\otimes\ldots\otimes v_{n+1}^*,
\qquad
t \, \mapsto\,
 ft\otimes v_1^*\otimes\ldots\otimes v_{n+1}^*,\ 
\end{equation}
$$
t^b\, \mapsto \, ft^b\otimes v_1^*\otimes\ldots\otimes v_{n+1}^*-
\sum_{l=0}^{b-2}\ \Big[et^l\otimes v_1^*\otimes\ldots\otimes 
v_n^*\otimes 2\sum_{i+j=b-l,\, i\geq j\geq 1}\ 
\left(\frac{e}{T^i} \frac{e}{T^j}v_{n+1}\right)^*+
$$
$$
ht^{l+1}\otimes v_1^*\otimes\ldots\otimes v_n^*\otimes 
\left(\frac{e}{T^{b-l-1}}v_{n+1}\right)^* \Big], 
$$
for $b\geq 2$.

\subsection{Theorem} 
\label{morphism} 
{\em Formulas (\ref{eta 1f})-(\ref{eta 0i}) define a monomorphism of complexes 
$$
\eta:\ \Omega^{\bullet}(U(z))[1]\lra 
C_{\bullet}(\sll_2(U(z));\otimes_{i=1}^{n+1}\ 
V^*_i)  .
$$}

\subsection{Beginning of the proof of Theorem \ref{morphism}} We should check that 
\begin{equation}
\label{equal}
\eta^1(\tilde d(x))=\mu(z)(\eta^0(x))
\end{equation} 
for any $x\in\Omega^0(U(z))$. We have 
\begin{equation}
1\, \overset{\eta}{\mapsto}\, f\otimes v_1^*\otimes\ldots
\otimes v_{n+1}^*\, \overset{\mu}{\mapsto}\,
\sum_{i=1}^n\ M^iv_1^*\otimes\ldots\otimes (fv_i)^*
\otimes\ldots\otimes v_{n+1}^*,
\end{equation}
$$
1\, \overset{\tilde d}{\mapsto}\, -\frac{1}{\kappa}\sum_{i=1}^n\ 
M^i\frac{d(t-z_i)}{t-z_i}\, \overset{\eta}{\mapsto}\,
\sum_{i=1}^n\ M^iv_1^*\otimes\ldots\otimes (fv_i)^*
\otimes\ldots\otimes v_{n+1}^*, 
$$
and these formulas agree with (\ref{equal}).
Next we have
\begin{equation}
t\, \overset{\eta}{\mapsto}\, ft\otimes v_1^*\otimes\ldots 
\otimes v_{n+1}^*\, \overset{\mu}{\mapsto}
\end{equation}
$$
\sum_{i=1}^n\ M^iz_i\, v_1^*\otimes\ldots\otimes (fv_i)^*
\otimes\ldots\otimes v_{n+1}^*\,+\,
(k-M^{n+1})v_1^*\otimes\ldots\otimes \left(\frac{e}{T}v_{n+1}\right)^*,
$$
$$
t\, \overset{\tilde d}{\mapsto}\, \frac{1}{\kappa}\Big[\left(\kappa-\sum_{i=1}
^n\ M^i\right)dt - \sum_{i=1}^n\ M^iz_i\frac{d(t-z_i)}{t-z_i}\Big]
\overset{\mu}{\mapsto} 
$$
$$
\left(k-\left(\sum_{i=1}^n\ M^i -2\right)\right) v_1^*\otimes\ldots\otimes 
\left(\frac{e}{T}v_{n+1}\right)^*\,+\,\sum_{i=1}^n M^iz_i  \,v_1^*\otimes
\ldots\otimes (fv_i)^*\otimes\ldots\otimes v_{n+1}^* ,
$$
and these formulas also agree with (\ref{equal}).
Notice that calculating the action on $M_{n+1}^*$ we use the automorphism $\pi$,
see formula (\ref{action}).

Similarly, to prove (\ref{equal}) for $x=t^2$, one needs 
the identity 
\begin{equation}
\label{b=2}
\left(2\kappa-\sum_{i=1}^n\ M^i\right)\left(\frac{e}{T^2}v_{n+1}\right)^*=
\frac{e}{T^2}\cdot v_{n+1}^*
-2f\left(\left(\frac{e}{T}\right)^2v_{n+1}\right)^*+\frac{h}{T}\left(\frac{e}{T}
v_{n+1}\right)^*,
\end{equation}
and to prove (\ref{equal}) for $x=\frac{1}{t-z_i}$, 
one needs the identity
\begin{equation}
\label{Ex}
(M^i+\kappa)\left(\frac{f}{T}v_i\right)^*=\frac{f}{T}\cdot v_i^*-
\frac{h}{T}(fv_i)^*-2\frac{e}{T}(f^2v_i)^* .
\end{equation}
\subsection{Theorem} 
\label{identities} {\em For $M,\ k\in\BC$, the following 
identities hold in the contragradient Verma module 
$V(M,k-M)^*$:\ 

{\em (a)}   for $b\geq 1$,  we have
$$
(M+b(k+2))\left(\frac{f}{T^b}v\right)^*=\frac{f}{T^b}\cdot v^*-
\sum_{l=1}^b\ \Big[2\frac{e}{T^l}\sum_{i+j=b-l,\ i\geq j\geq 0
\atop }\ 
\left(\frac{f}{T^i}\frac{f}{T^j}v\right)^*+\frac{h}{T^l}\left(
\frac{f}{T^{b-l}}v\right)^*\Big];
$$

{\em (b)}  for $b\geq 2$, we have
$$
(k-M+(b-1)(k+2))\left(\frac{e}{T^b}v\right)^*=\frac{e}{T^b}\cdot v^*+
\sum_{l=0}^{b-2}\ \Big[-2\frac{f}{T^l}\sum_{i+j=b-l,\ i\geq j
\geq 1}\ \left(\frac{e}{T^i}\frac{e}{T^j}v\right)^*+\frac{h}{T^{l+1}}
\left(\frac{e}{T^{b-l-1}}v\right)^*\Big].
$$}

\subsection{Proof of the theorem}

The theorem is proved by direct verification. Each term of the  expression in (a) is of degree
$(b+1,b)$. The basis of $V(M,k-M)_{(b+1,b)}$ is described in Section \ref{BASIS}. This gives us the dual basis
of  $V(M,k-M)_{(b+1,b)}^*$.  One  calculates in straightforward way
the right-hand side of (a) in that basis and obtains the left-hand side.
For example, for $b=1$, the space  $V(M,k-M)_{(2,1)}^*$ has the basis
$(\frac fTv)^*$, $(f \frac hTv)^*$,  $(f^2\frac eTv)^*$, and we have
$\frac fT\cdot v^* = (M+k)(\frac fTv)^* + (2M-2k)(f^2\frac eTv)^*$,
$\frac hT\cdot(fv)^*= -2(\frac fTv)^* + 4(f \frac hTv)^* + (4M-4k)(f^2\frac eTv)^*$,
$\frac eT(f^2v_i)^* = -2(f \frac hTv)^* - (2M-2k) (f^2\frac eTv)^*$. 
 By adding these expressions we get the formula
$(M+\kappa)\left(\frac{f}{T}v\right)^*=\frac{f}{T}\cdot v^*-
\frac{h}{T}(fv)^*-2\frac{e}{T}(f^2v)^*$ which gives statement (a) for $b=1$ and formula (\ref{Ex}).

Similarly each term of the expression in (b) is of degree
$(b-1,b)$. The basis of $V(M,k-M)_{(b-1,b)}$ is described in Section \ref{BASIS}. This gives us the dual basis
of  $V(M,k-M)_{(b-1,b)}^*$.  One  calculates in straightforward way
the right-hand side of (b) in that basis and obtains the left-hand side.
For example, for $b=2$ the space  $V(M,k-M)_{(1,2)}^*$ has the basis
$((\frac e T)^2f v)^*, (\frac eT\frac hT v)^*, (\frac e{T^2}v)^*$ and we have
$\frac e{T^2}\cdot v^*=-2M((\frac e T)^2f v)^*-2k (\frac eT\frac hT v)^*+(2k-M)(\frac e{T^2}v)^*$,
$f\cdot((\frac eT)^2v)^* = M((\frac e T)^2f v)^*-2 (\frac eT\frac hT v)^*$,
$\frac hT\cdot(\frac eTv)^* =
4M((\frac e T)^2f v)^*+(2k-4) (\frac eT\frac hT v)^*+2(\frac e{T^2}v)^*$.
 By adding these expressions  we get the formula
$(2k+2-M)(\frac e{T^2}v)^* = \frac e{T^2}\cdot v^*-2f\cdot((\frac eT)^2v)^*+
\frac hT\cdot(\frac eTv)^*$ which gives statement (b) for $b=2$ and  formula 
(\ref{b=2}).

The complete proofs of the theorem will be published elsewhere.
$\Box$

\subsection{End of the proof of Theorem \ref{morphism} } Theorem \ref{morphism} is a direct corollary 
of Theorem \ref{identities}, cf. (\ref{bff}), 
(\ref{bfi}) and (\ref{eta 1f})-(\ref{eta 0i}).      $\Box$

\section{Singular vectors in Verma modules} 

\subsection{} Let $S:\ V(M,k-M)\lra V(M,k-M)^*$ be the Shapovalov form. 
Set
\begin{equation}
X_b(M,k-M):= S^{-1}\left((M+b(k+2))\left(\frac{f}{T^b}v\right)^*\right),\
\end{equation}
$$ 
Y_b(M,k-M):= S^{-1}\left((k-M+(b-1)(k+2))\left(\frac{e}{T^b}v\right)^*\right).
$$
For generic values of $M$ and $k$, the Shapovalov form $S$ is non-degenerate 
and $X_b$ and $Y_b$ are well defined elements of the Verma module 
$V(M,k-M)$. 
The basis in $V(M,k-M)$ allows us to compare these vectors for different 
values of $k, M$. Obviously, $X_b(M,k-M),\ Y_b(M,k-M)$ are 
holomorphic functions of $k, M$ for generic $k, M$. 

Consider the resonance lines
\begin{equation}
M=l-1-(a-1)(k+2),
\qquad M=-l-1+a(k+2),\qquad k+2=0 ,
\end{equation}
$(l, a\in\BZ_{>0}$) on the $(M,k)$-plane, cf. \ref{reduc}. 

\subsection{Theorem} 
\label{limit}
\begin{enumerate}
\item[(a)]
 {\em  Let $b\geq 0$ and let
$(M_0, k_0)$ be a point
of the line $\{ M=-b(k+2) \}$ which does not belong to  other
resonance lines. Then the vector-valued
function $X_b(M,k-M)$ can be analytically
continued to the point $(M_0,k_0)$ and the vector $X_b(M_0,k_0-M_0)$ is a
(nonzero) singular vector of $V(M,k-M)$.

\item[(b)] 
Let $b>0$ and let $(M_0,k_0)$ be a point of the line
$\{ M=-2+b(k+2) \}$ which does not belong to  other resonance
lines. Then the vector-valued
function $Y_b(M,k-M)$ can be analytically
continued to the point $(M_0,k_0)$ and the vector $Y_b(M_0,k_0-M_0)$ is a
(nonzero) singular vector of $V(M,k-M)$.}
                                        
\end{enumerate}

{\it Proof} of (a). According to Theorem \ref{identities}, 
\begin{equation}
\label{aux}
X_b(M,k-M)=\frac{f}{T^b}v-\sum_{l=1}^b\ \Big[2\frac{e}{T^l}
\sum_{i+j=b-l,\ i\geq j\geq 0}\ S^{-1}\left(\frac{f}{T^i}\frac{f}{T^j}v\right)^*+
\frac{h}{T^l}S^{-1}\left(\frac{f}{T^{b-l}}v\right)^*\Big].
\end{equation}
The right-hand side can be analytically continued to $(M_0,k_0)$
since the elements
\\
 $S^{-1}\left(\frac{f}{T^i}\frac{f}{T^j}v\right)^*$
and $S^{-1}\left(\frac{f}{T^{b-l}}v\right)^*$ are well defined at  
$(M_0,k_0)$ by the results in Section \ref{reduc}. 
We have 
$S(X_b(M_0,k_0-M_0))=0$ by definition. Let us check that 
$X_b(M_0,k_0-M_0)\neq 0$. In fact, consider the basis 
$\{\frac{e}{T^{l_1}}\cdot\ldots\cdot\frac{e}{T^{l_{\alpha}}}
\frac{h}{T^{j_1}}\cdot\ldots\cdot\frac{h}{T^{j_{\beta}}}
\frac{f}{T^{i_1}}\cdot\ldots\cdot\frac{f}{T^{i_{\gamma}}}\}$ in 
$V(M_0,k_0-M_0)$. Formula (\ref{aux}) shows that the basis vector 
$\frac{f}{T^b}v$ comes to $X_b$ with coefficient $1$. This proves 
part (a) of the theorem. Part (b) is proved similarly. 
\qed

\subsection{} 
\label{MFF vectors} 
{\bf Corollary of Theorem \ref{limit}.}
{\em If $b\geq 0$ and 
$(M_0, k_0)$ is a point
of the line $\{ M=-b(k+2) \}$, which does not belong to  other
resonance lines, then   $X_b(M_0,k_0-M_0)$
is proportional to the Malikov-Feigin-Fuchs vectore $F_{12}(1,b+1, k_0+2)$.
Similarly, if $b>0$ and  $(M_0,k_0)$ is a point of the line
$\{ M=-2+b(k+2) \}$ which does not belong to  other resonance
lines, then $Y_b(M_0,k_0-M_0)$ is proportional to 
the Malikov-Feigin-Fuchs vectore $F_{21}(1,b, k_0+2)$.
$\Box$

}

\subsection{} 
\label{res sing} 
{\bf Corollary of formulas  (\ref{bff}) and (\ref{bfi}).}
{\em If the resonance condition $k-M^{n+1}+(b-1)\kappa=0$ of the identity
\ref{identities}.b holds, then formula  (\ref{bfi}) gives the following 
cohomological relation
between the logarithmic forms $\omega_j$ introduced in \ref{log forms}:
\begin{equation}
\sum_{m\geq 0}\, 
\sum_{l_0+\ldots +l_m=b,
\atop
l_0,\ldots,l_m >0}
\, (-\kappa)^{-m}\,
(\,\prod_{i=1}^m \,(\, \frac{1}{l_i}\, \sum_{j=1}^n\, z_j^{l_i}\,M^j\,)\,)\,
\sum_{j=1}^n\, z_j^{l_0}\,\omega_j\,\sim \, 0\, ,
\end{equation}
where $\kappa=k+2$. Similarly if the resonance conditiion $M^p+b\kappa=0$
 of the identity \ref{identities}.a holds for some $p\leq n$,  then formula  (\ref{bff}) 
induces the following  cohomological relation
between the logarithmic forms $\omega_j$:
$$
\sum_{m\geq 0}\, [\,
\sum_{l_0+\ldots +l_m=b,\atop
l_0,\ldots, l_m>0}
\, (-\kappa)^{-m} \,
(\,\prod_{i=1}^m \,(\frac{1}{l_i}\,
\sum_{
j=1,\ldots, n, \atop j\neq p}
\, \frac{M^j}{(z_j-z_p)^{l_i}}\,)\,)\, 
\sum_{j=1,\ldots, n,\atop j\neq p}\,
\, \frac{\omega_j}{(z_j-z_p)^{l_0}} -
$$
\begin{equation}
\sum_{l_1+\ldots + l_m=b,\atop l_1,\ldots, l_m>0 }
\, 
(-\kappa)^{-m}(\prod_{i=1}^m\,(\frac{1}{l_i}
\sum_{j=1,\ldots n, \atop j\neq p }
\, \frac{M^j}{(z_j-z_p)^{l_i}}
))\omega_p ]\sim 0 ,
\end{equation} 
cf. Examples \ref{e1}, \ref{e2}. 
$\Box$
}

For instance, if $M^1+\kappa=0$, then 
\begin{equation}
\sum_{j>1}\ \frac{\omega_j}{z_j-z_1} + 
\frac{1}{\kappa} \left(\sum_{j>1}\ \frac{M^j}{z_j-z_1}\right)\omega_1 
\sim 0\ .
\end{equation}

\end{document}